\documentclass{article}
\usepackage{amssymb,latexsym,amsmath,amsthm,mathrsfs}
\usepackage{graphicx}
\newcommand{\comment}[1]{}

\begin{document}
\title{Observations on two fourth powers whose sum is equal to the sum of two
other fourth powers\footnote{Presented to the St. Petersburg Academy
on January 13, 1772. Originally published as
{\em Observationes circa bina biquadrata quorum summam in duo alia
biquadrata resolvere liceat},
Novi Commentarii Academiae scientiarum Imperialis Petropolitanae
\textbf{17} (1773), 64--69.
E428 in the Enestr{\"o}m index.
Translated from the Latin by Jordan Bell,
Department of Mathematics, University of Toronto, Toronto, Ontario, Canada.
Email: jordan.bell@gmail.com}}
\author{Leonhard Euler}
\date{}
\maketitle

\begin{center}
{\Large Summarium}
\end{center}
{\small Among the theorems concerned with the properties of numbers, one often
tries to prove that the sum of three fourth powers cannot also be a fourth power,
that is, that the equation
\[
a^4+b^4+c^4=d^4
\]
is impossible; but the same does not also hold for differences, since the
equation
\[
A^4+B^4-C^4=D^4
\]
or
\[
A^4+B^4=C^4+D^4
\]
or
\[
A^4-D^4=C^4-B^4
\]
can be resolved even in infinitely many ways; and even if perhaps this
has already been found by other Geometers, the method which the Celebrated
Author used here has been found for more than reason
to merit the full attention of Analysts. The four smallest numbers satisfying the proposed
question were found thus by this Illustrious Man:\footnote{Translator: As I mention in the footnote at the end of
\S 8, there was a mistake in the calculations in \S 8 and these numbers do
not
satisfy the equation.}
\[
A=477069, \quad B=8497, \quad C=310319 \quad \textrm{and} \quad
D=428397.
\]
}

1. Since it has been demonstrated that neither the sum nor the difference
of two fourth powers can be a square, much less can
it be a fourth power;\footnote{Translator: See
Andr\'e Weil, {\em Number theory: an approach through history}, Birkh\"auser,
1984,
p. 79.}
one would hardly have any more doubt that the sum of three fourth powers
could ever be a fourth power, even if this has still not been
demonstrated.\footnote{Translator: This statement turns out to be false.
Noam Elkies in {\em On $A^4+B^4+C^4=D^4$}, Math. Comp. \textbf{51} (1988),
pp. 825--835, finds the counterexample
$2682440^4+15365639^4+18796760^4=20615673^4$. See also
Lawrence D'Antonio, {\em Teaching elliptic curves using original sources},
From calculus to computers: using the last 200 years of mathematics history
in the classroom (Amy Shell-Gellasch and Dick Jardine, eds.), Mathematical
Association of America, 2005, pp. 25--44.}
We are rightly uncertain whether one can find four fourth powers whose sum is
a fourth power, since so far no one
has exhibited such fourth powers.\footnote{Translator: Euler is unable
to find four fourth powers whose sum is a fourth power, but believes that
it is indeed possible. See Euler's August 4, 1753 letter to Goldbach.}

2. Even if it could be demonstrated that one cannot give three fourth powers
whose sum was a fourth power, this could by no means be extended
to differences: it cannot be therefore affirmed that an equation such as
\[
A^4+B^4-C^4=D^4
\]
is impossible; for I have observed that this equation can be
resolved in infinitely many ways.
However I do not wish to claim that nobody has hitherto shown this,
but now at least there is not the opportunity to unravel
all the annals of this type of Analysis; whatever the case, I believe
the method
which I have used 
to be not unworthy of the utmost attention. It is apparent that this question revolves
around two pairs of fourth powers, whose sums or differences are equal;
for if it were the case that
$A^4+B^4=C^4+D^4$ then it will also be
the case that $A^4-D^4=C^4-B^4$. This suggested to us the following Problem.

\begin{center}
{\Large Problem}
\end{center}
{\em To find two fourth powers $A^4$ and $B^4$ whose sum may be
resolved into two other fourth powers, so that the equality}
\[
A^4+B^4=C^4+D^4
\]
{\em is obtained.}

\begin{center}
{\Large Solution\footnote{Translator: Thomas L. Heath, {\em Diophantus of Alexandria: a study in
the history of Greek algebra}, second ed., Cambridge, 1910,
problem 17, pp. 377--380 gives an exposition of Euler's solution.}}
\end{center}

3. Since here it ought to be that
\[
A^4-D^4=C^4-B^4,
\]
let us put
\[
A=p+q,\quad D=p-q,\quad C=r+s \quad \textrm{and} \quad B=r-s,
\]
so that the following nicer equation is produced
\[
pq(pp+qq)=rs(rr+ss),
\]
which indeed is evidently satisfied by taking $r=p$
and $s=q$; truly we gained nothing from this since the case that arises
$C=A$ and $B=D$ is obvious by itself; but at the same time,
this very case avails to lead to other solutions.\footnote{Translator: Of course
the specific case $r=p$ and $s=q$ only yields one solution. I presume
this means that by trying to modify this case by making
$r$ and $p$ multiples of each other and $s$ and $q$ multiples of each other we are led to other solutions.}

4. Now let us set
\[
p=ax, \quad q=by, \quad r=kx \quad \textrm{and} \quad s=y,
\]
so that the following equation is obtained which needs to be resolved
\[
ab(aaxx+bbyy)=k(kkxx+yy),
\]
from which we deduce at once
\[
\frac{yy}{xx}=\frac{k^3-a^3b}{ab^3-k},
\]
which fraction has therefore to be made a square.\footnote{Translator:
A square of a rational number, not necessarily a square of an integer.}
Here the case immediately presents itself which occurs by taking
$k=ab$; for then it would be
\[
\frac{yy}{xx}=\frac{a^3b(bb-1)}{ab(bb-1)}=aa,
\]
whence it would happen that $y=a,x=1$,\footnote{Translator: 
$x$ and $y$ are assumed to be relatively prime, for otherwise
$p,q,r,s$ would have a common factor, and thus $A,B,C,D$ would have a common
factor. But $A,B,C,D$ can be assumed not to have a common factor.}
and then $p=a,q=ab,r=ab,s=a$, and these values only yield the obvious
case.

5. However, further pursuing this case
let us set $k=ab(1+z)$, and our equation will be transformed into this
form
\[
\frac{yy}{zz}=\frac{a^3b(bb-1+3bbz+3bbz^2+bbz^3)}{ab(bb-1-z)}=
aa\frac{bb-1+3bbz+3bbz^2+bbz^3}{bb-1-z}
\]
and from this equation we elicit
\[
\frac{y}{x}=\frac{a\surd((bb-1)^2+(bb-1)(3bb-1)z+3bb(bb-2)zz+bb(bb-4)z^3-bbz^4)}{bb-1-z}.
\]
Then to make the formula
\[
(bb-1)^2+(bb-1)(3bb-1)z+3bb(bb-2)zz+bb(bb-4)z^3-bbz^4
\]
a square, let us put its square root
\[
=bb-1+fz+gzz
\]
and we shall take the letters $f$ and $g$ such that the first three terms
cancel; so since the square of this form is
\[
(bb-1)^2+2(bb-1)fz+2(bb-1)gzz+ffzz+2fgz^3+ggz^4,
\]
the first terms spontaneously 
cancel; so that the same happens for the second terms, one ought to take
\[
f=\frac{3bb-1}{2}
\]
and for the third we will have
\[
3bb(bb-2)=2(bb-1)g+\frac{9b^4-6bb+1}{4},
\]
from which one gets
\[
g=\frac{3b^4-18bb-1}{8(bb-1)},
\]
and with these values settled the equation to be solved is
\[
(gg+bb)z=bb(bb-4)-2fg,
\]
from which we get
\[
z=\frac{bb(bb-4)-2fg}{bb+gg}.
\]

6. Then the letter $b$ still remains our choice; it can therefore be assumed
as we please and we will have simultaneously determined the quantity $z$;
then
we will immediately have
\[
x=bb-1-z \quad \textrm{and} \quad y=a(bb-1+fz+gzz)
\]
and then in turn
\[
\begin{split}
&p=a(bb-1-z), \quad q=ab(bb-1+fz+gzz),\\
&r=ab(1+z)(bb-1-z), \quad s=a(bb-1+fz+gzz);
\end{split}
\]
since all these formulas are divisible by $a$ it can be eliminated by dividing,
so that it will be
\[
\begin{split}
&p=bb-1-z,\quad q=b(bb-1+fz+gzz),\\
&r=b(1+z)(bb-1-z),\quad s=bb-1+fz+gzz,
\end{split}
\]
where it should be noted that if the numbers $x$ and $y$ had a common factor,
it could eliminated before the letters $p,q,r,s$ are determined.

It will thus be worthwhile to work out some special solutions;
and indeed it is at once apparent that one cannot take $b=1$, because this would
make $g=\infty$; much less can one put $b=0$, because this would make
$q=0$. Therefore we shall show the two cases $b=2$ and $b=3$.

\begin{center}
{\Large First special solution}
\end{center}

7. Let $b=2$, and the above values are gathered as follows
\[
f=\frac{11}{2},\quad g=-\frac{25}{24}, \quad z=\frac{6600}{2929};
\]
then since the letter $a$ does not enter into the calculation, unity
can be written in its place; then indeed it will be
\[
\begin{split}
&x=3-\frac{6600}{2929}=\frac{2187}{2929},\\
&y=3+\frac{11}{2}\cdot \frac{6600}{2929}-\frac{25}{24}\cdot \frac{6600^2}{2929^2}
=3+\frac{55407\cdot 1100}{2929^2}=\frac{3\cdot 28894941}{2929^2}.
\end{split}
\]
Thus the whole matter reduces to the relation between $x$ and $y$.
Since 
\[
\frac{y}{x}=\frac{3\cdot 28894941}{2187\cdot 2929}=\frac{28894941}{2929\cdot 729}
=\frac{3210549}{2929\cdot 81}
=\frac{1070183}{27\cdot 2929},
\]
we will have
\[
x=79083 \quad \textrm{and} \quad y=1070183;
\]
then because
\[
k=2(1+z)=\frac{2\cdot 9529}{2929}=\frac{19058}{2929}
\]
we conclude that
\[
\begin{split}
&p=79083,\quad q=2\cdot 1070183=2140366,\\
&r=27\cdot 19058=514566,\quad s=1070183.
\end{split}
\]
Consequently we obtain for the roots of the fourth powers
\begin{eqnarray*}
A=p+q=2219449,&C=r+s=1584749,\\
B=r-s=-555617,&D=p-q=-2061283
\end{eqnarray*}
and then $A^4+B^4=C^4+D^4$.

\begin{center}
{\Large Second special solution}
\end{center}

8. Let $b=3$, and it will be
\[
f=13, \quad g=\frac{5}{4}, \quad \textrm{and then} \quad z=\frac{200}{169}
\] 
and hence
\[
k=\frac{3\cdot 369}{169}=\frac{1107}{169}=\frac{9\cdot 123}{169}
=\frac{27\cdot 41}{169};
\]
further,
\[
x=\frac{8\cdot 144}{169}=\frac{128\cdot 9}{169}
\]
and\footnote{Translator: The {\em Opera omnia} corrects an error
in the original edition of the paper, which had $y=\frac{8\cdot 150911}{169^2}$.}
\[
y=8+\frac{200}{169}\Big(13+\frac{5}{4}\cdot \frac{200}{169}\Big)
=8+\frac{200}{169}\cdot \frac{2447}{169}=\frac{8\cdot 89736}{169^2}
\]
and thus it will be
\[
x:y=8\cdot 144\cdot 169:8\cdot 89736=6\cdot 169:3739
\]
and so
\[
x=6\cdot 169=1014 \quad \textrm{and} \quad y=3739,
\]
from which values we get
\begin{eqnarray*}
p=1104,&r=6642=6\cdot 1107,\\
q=11217,&s=3739.
\end{eqnarray*}
And then the letters $A,B,C,D$ are themselves collected
\begin{eqnarray*}
A=12231,&C=10381,\\
B=2903,&D=-10203
\end{eqnarray*}
and it will again be
\[
A^4+B^4=C^4+D^4
\]
and these seem to be the smallest numbers satisfying our question.\footnote{Translator: Euler made this statement about the numbers he got from his mistaken
calculation, $A=477069,B=8497,C=310319,D=42897$. The {\em Opera omnia} notes
that Euler corrects this error later in his
1780
{\em Dilucidationes circa binas summas duorum biquadratorum inter se aequales},
E776, where he gives the solution $A=542, B=103,
C=359, D=514$. For a simpler parametric solution that yields the smaller
solution $A=133,B=134,C=158,D=59$, see Hardy and Wright,
{\em An introduction to the theory of numbers}, fifth ed., Oxford, p. 201.}

\end{document}